%% file: main.tex
\documentclass[12pt,reqno]{amsart}
\usepackage{fullpage}
\usepackage{url}
\usepackage{xcolor}

\include{pack}

\begin{document}
\title{Towards the weighted bounded negativity conjecture for blow-ups of algebraic surfaces}
\author{Roberto Laface \and Piotr Pokora}
\date{\today}


\address{Roberto Laface \newline Technische Universit\"at M\"unchen, Zentrum Mathematik - M11,
    Boltzmannstra{\ss}e 3,
    D-85478 Garching bei M\"unchen, Germany
}
    \email{laface@ma.tum.de}

\address{Piotr Pokora \newline
   Institut f\"ur Algebraische Geometrie,
    Leibniz Universit\"at Hannover,
    Welfengarten 1,
    D-30167 Hannover, Germany. \\
Current address: Institute of Mathematics,
   Polish Academy of Sciences,
   ul. \'{S}niadeckich 8,
   00-656 Warszawa, Poland.}
   \email{piotrpkr@gmail.com, ppokora@impan.pl}

\maketitle

\thispagestyle{empty}
\begin{abstract}
In the present paper we focus on a weighted version of the Bounded Negativity Conjecture, which predicts that for every smooth projective surface in characteristic zero the self-intersection numbers of reduced and irreducible curves are bounded from below by a function depending on the intesection of curve with an arbitrary big and nef line bundle that is positive on the curve. We gather evidence for this conjecture by showing various bounds on the self-intersection number of curves in an algebraic surface. We focus our attention on blow-ups of algebraic surfaces, which have so far been neglected.
\end{abstract}

\section*{Introduction}
In the last years, negative curves on surfaces have been researched extensively because of their connection to many open problems. Among these, one cannot refrain from mentioning Nagata's conjecture \cite{Nagata} or the SHGH conjecture \cite{Ciro}. The present paper is devoted to yet another open question in the geometry of complex surfaces:

\begin{conj}[Bounded Negativity Conjecture]\label{bnc} For every smooth projective surface $X$ over the complex numbers, there exists a nonnegative integer $b(X) \in \Z$ such that $C^2 \geq -b(X)$ for all integral curves $C \subset X$.
\end{conj}

The Bounded Negativity Conjecture (BNC in short) has a long oral tradition, and it seems to date back to F. Enriques. In some cases, the conjecture is known to hold true, for instance when the anti-canonical bundle is $\Q$-effective or when the surface is equipped with a surjective endomorphism of degree $d >1$. However, if one considers non-minimal surfaces, e.g.~blow-ups of a surface for which BNC is known to hold, then very little is known and the problem acquires a very different flavor.

As it turns out, the BNC is equivalent to the statement in Conjecture \ref{bnc} where one allows $C$ to be any reduced curve in $X$ \cite[Proposition 3.8.2]{RDLS}. This has paved the way to the study of the BNC from the point of view of configurations of curves via the notion of H-constant \cite{BdRHHLPSz}. The H-constant is an asymptotic invariant that has the potential of studying the BNC on all blow-ups of a given algebraic surface at all possible configurations of points on it simultaneously, see for instance \cite{BdRHHLPSz, LP2015,LP2016, Roulleau, Roulleau1, PR}.

In the present paper, we go back to focusing our attention on integral curves and bounding their negativity. In \cite[Conjecture 3.7.1]{RDLS}, the authors formulated the following variant of the BNC.

\begin{conj}[Weighted BNC]\label{wbnc} For every smooth projective surface $X$ over the complex numbers, there exists a nonnegative integer $b_w \in \Z$ such that $C^2 \geq - b_w(X) \cdot (C.H)^2$ for all integral curves $C \subset X$ and all big and nef line bundles $H$ for which $C.H > 0$.	
\end{conj}

Notice that we are not asking for the self-intersection of a curve $C$ to be bounded from below, but rather that the \textit{weighted self-intersection} $C^2 /(C.H)^2$ of $C$ be so, hence the adjective "weighted". Put differently, Conjecture \ref{wbnc} is asking for a bound on the self-intersection of all integral curves on $X$ that depends on both $X$ and the degree of the curve $C$ with respect of every big and nef line bundle over which the curve is positive.  The importance of the weighted BNC lies in the fact that it implies positivity of the global Seshadri constant of ample line bundles at all points of a given surface $X$ \cite[Proposition 3.6.2]{RDLS}. 

Our paper aims at gathering evidence for the validity of this conjecture. More precisely, we provide bounds for the self-intersection numbers of irreducible and reduced curves on blow-ups of algebraic surfaces at mutually distinct points. The bounds depend on the degree of the curve with respect to an explicitly constructed big and nef line bundle $\Gamma$, and in fact it holds for the cone $\text{Nef}(X)+\Gamma$ (the translate of the nef cone by $\Gamma$). 

The technical heart is Theorem \ref{linebundle}, where we construct a line bundle on a blow-up $Y$ of $X$ at $n$ distinct points that naturally arises from $X$. We prove this result by first showing a generalization of a result due to Sakai \cite{Sakai} and Orevkov-Zaidenberg \cite{OZ95}, together with estimates on the Milnor numbers of isolated singularities. This provides a function that depends linearly on the degree with respect to a given line bundle, while the conjecture only predicts that such a function should be quadratic. 

Our results give a uniform treatment of the case of surfaces of non-negative Kodaira dimension (see Corollary \ref{thmA}):
\begin{theoremA}
		Assume $X$ is a surface of non-negative Kodaira dimension and let $f: Y\rightarrow X$ be the blowing up of $X$ along $n$ mutually distinct points. Then there exists a big and nef line bundle $\Gamma$ that bounds negativity linearly, i.e.,
        \[ C^2 \geq -\frac{1}{2} \big( \delta(X) + C.\Gamma \big) - n, \]
        for every integral curve $C \subset Y$, where $\delta(X)=3e(X)-K_{X}^{2}$.
\end{theoremA}

Turning to surfaces of Kodaira dimension $\kappa=-\infty$, we are able to give a very neat picture in the case of blow-ups of $\IP^2$ (see Theorem \ref{boundrat}):

\begin{theoremB}
		Let $\sigma: Y \longra \IP^2$ be the blow-up of $\IP^2$ at $n$ mutually distinct points in $\IP^2$, and let $C$ be an irreducible and reduced curve on $Y$. Then,
        \[ C^2 \geq - 2n(C.L), \]
        where $L$ is the pull-back of a line in $\IP^2$.
\end{theoremB}

We also have partial results on blow-ups of Hirzebruch surfaces, and we refer to Section \ref{hirzebruch} for the details. We are working exclusively over the complex numbers.

\section{Generalization of a result of Sakai and Orevkov-Zaidenberg}

In this section, we are going to provide a generalization of the following result, proven independently by Sakai \cite{Sakai} and Orevkov-Zaidenberg \cite{OZ95}.

\begin{thm}\label{SOZ}
Let $C$ be a reduced and irreducible curve in $\mathbb{P}^{2}$ of degree $d$ having singular points $p_{1}, ...,p_{s}$. We denote by $m_{p_{i}}$ and $\mu_{p_{i}}$ the corresponding multiplicity and the Milnor number of $p_{i}$. If the logarithmic Kodaira dimension of $\mathbb{P}^{2} \setminus C$ is non-negative, then
$$\sum_{i=1}^{s} \bigg(1 + \frac{1}{2m_{p_{i}}}\bigg)\mu_{p_{i}} \leq d^{2} - \frac{3}{2}d.$$
\end{thm}

For the definition of the Milnor number of a singularity we refer to \cite[\S 7]{Milnor}. Our aim is to show that the above inequality holds true in a broader setting. Before we present the result, let us recall that one has the following variation on Max Noether's inequality \cite[Satz~5, p.~835]{KB}.

\begin{thm}\label{sakaigeneral}
Let $X$ be a smooth complex projective surface and $C \subset X$ an irreducible and reduced curve with singular points $p_{1},...,p_{s}$ and denote by $\mu_{p}$ the Milnor number of the singularity of $C$ at $p$. If $K_{X}$ denotes the canonical divisor of $X$ and $e(C)$ denotes the topological Euler characteristic of $C$, then
$$e(C) = \sum_{p \in Sing(C)} \mu_{p} -(C.C + C.K_{X}).$$
\end{thm}

Our approach is to follow an idea of Sakai \cite[\S 1]{Sakai}, which we illustrate below. 

\begin{constr}\label{constr}
For an irreducible and reduced curve $C$ on a smooth surface $X$, we denote by $f: S \rightarrow X$ the minimal sequence of blow-ups such that the (reduced) total transform of $C$ has normal crossings. Let $\{E_{1}, ..., E_{n}\}$ be the set of exceptional curves for $f$ (i.e.~the exceptional divisors that arise when performing the blow-up $f$), and we set $D = \tilde{C} + \sum_{i} E_{i}$, $\tilde{C}$ being the strict transform of $C$ along $f$. For a singularity $(C,p)$:
\begin{enumerate}
\item $m_{p}$ is the multiplicity of $(C,p)$;
\item $r_{p}$ is the number of branches of $(C,p)$;
\item $E$ is the reduced exceptional divisor of $f$, $E = \sum_i E_i = \sum_{p \in Sing(C)} E_p$, where $E_p$ is the reduced exceptional divisor over the point $p \in X$;
\item $\omega_{p} = -E_{p}^{2}$. 
\end{enumerate}

Let us recall that $f^*C = \tilde{C} + \sum_i m_{i} \bar{E}_{i}$, $\bar{E}_{i}$ being the total transform of $E_i$ in $S$, while the reduced exceptional divisor satisfies $E.\tilde{C}= \sum_p r_p$. Indeed, let us consider a singularity $ (C,p) $: along a resolution, the $ r_p $ branches will get separated. As blowing-up is an isomorphism outside of the center, for every branch there exists unique exceptional divisor intersecting it transversally. The intersection point of the exceptional divisor and the branch maps to $ p $. By summing over all singular points, we obtain the desired formula.
\end{constr}

\begin{defi}
For a singularity $(C,p)$ we denote by $(m_{1} = m_{p}, m_{2}, ..., m_{n})$ the sequence of multiplicities of all infinitely near points of $p$ in $f$. We set
$$\eta_{p} = \sum_{j=1}^{n} (m_{j}-1),$$
and since $\sum_{j}m_{j}(m_{j}-1) = \mu_{p} + r_{p} - 1$ by \cite[p.85]{Milnor}, then we have
$$\sum_{j} (m_{j}^{2} - 1) = \mu_{p}+r_{p}-1 + \eta_{p}.$$
\end{defi}

We are now ready to show our version of the Orevkov-Sakai-Zaidenberg inequality, which we will employ in the study of the negativity of a surface carried out in Section \ref{bounding_negativity}.

\begin{thm}\label{SOZgen} Let $C$ be an irreducible and reduced curve in a smooth complex projective surface $X$ having singular points $p_{1},...,p_{s}$. We denote by $m_{p_i}$ and $\mu_{p_i}$ the corresponding multiplicities and the Milnor numbers of $p_{i}$'s. Assume that the logarithmic Kodaira dimension of $X\setminus C$ is non-negative, then one has
$$\sum_{p \in Sing(C)} \bigg( 2 + \frac{1}{m_{p}}\bigg) \mu_{{p}} \leq 3e(X) - K_{X}^{2} + 2C^{2} +K_{X}.C.$$
\end{thm}
\begin{proof}
Since $|m(K_{S} + D)| \neq \emptyset$ for a certain positive integer $m$, we can use the logarithmic Miyaoka-Sakai inequality \cite{Sakai1} for the pair $(S,D)$ as in Construction \ref{constr}, namely
$$(K_{S}+D)^{2} \leq 3(e(S)-e(D)).$$
First of all, we have
$$e(S)-e(D) = e(X) - e(C).$$
Now we would like to compute $(K_{S}+D)^{2}$. Following the idea of Sakai \cite[p.263]{Sakai}, we can see that:

\begin{align*}
D^2 &= (\tilde{C} + \sum_i E_i)^2 = \tilde{C}^2 - \sum_p (\omega_p -2r_p)\\
p_{a}(D) &= p_{a}(\tilde{C}) + \sum_p p_a(E_p) - s + \tilde{C}.E= p_{a}(\tilde{C}) +\sum_p (r_p-1)\\
e(\tilde{C}) &= e(C) + \sum_p (r_p -1) \\
(K_S+D).D &= 2p_{a}(D) -2 = 2p_{a}(\tilde{C}) -2 + 2 \sum_p (r_p - 1) = -e(C) + \sum_p(r_p -1)\\
K_S^2 - \tilde{C}^2 &= K_X^2 -C^2 + \sum_i (m_i^2-1).
\end{align*}
This leads to
$$ (K_S+D)^2 = K_X^2 -C^2 -2e(C) + \sum_p  (\mu_p + \omega_p + r_p -3 + \eta_p), $$
which implies
$$K_{X}^{2} - C^{2} - 2e(C) + \sum_p( \mu_p + \omega_p + r_p - 3 + \eta_{p} ) \leq 3e(X) - 3e(C),$$
by the logarithmic Miyaoka-Yau inequality. The above statement is equivalent to
$$e(C) +  \sum_p(\mu_p + \omega_p + r_p - 3 + \eta_{p}) \leq 3e(X) - K_{X}^{2} + C^{2}.$$
We have 
\begin{align*}
(K_{X}+C).C & = 2p_{a}(C)-2 = 2p_{a}(\tilde{C})-2 + \sum_{i} m_{i}(m_{i}-1)\\
&= -e(\tilde{C}) + \sum_{i} m_{i}(m_{i}-1) = -e(C) - \sum_p (r_p-1) +\sum_{i} m_{i}(m_{i}-1)\\
&= -e(C) - \sum_p \Big( r_p -1 - \sum_\text{$m_i$ over $p$} m_i(m_i-1) \Big) =-e(C) + \sum_p \mu_p,
\end{align*}
where the last equality follows from Milnor's formula \cite[p.85]{Milnor} and $\sum_\text{$m_i$ over $p$}$ means that we are summing up the multiplicities of the infinitely near points of $p$. From this, one has
$$\sum_p(2 \mu_p + \omega_p+ r_p - 3 + \eta_p) \leq 3e(X) - K_{X}^{2} + 2C^{2} + K_{X}.C.$$
As it was pointed out explicitly in \cite{OZ95}, we have the following inequality
$$\eta_p + \omega_p+r_p-3 \geq \mu_p/m_p.$$
This implies
$$\sum_{p \in Sing(C)} \bigg( 2 + \frac{1}{m_p}\bigg)\mu_p \leq 3e(X) - K_{X}^{2} + 2C^{2} + K_{X}.C,$$
which completes the proof.
\end{proof}


\section{Bounding negativity on surfaces with $\kappa \geq 0$}\label{bounding_negativity}

In this section, we would like to bound the negativity of curves on an algebraic surface, having in mind the Weighted BNC as a goal. Let $X$ a smooth projective surface over the complex numbers, and let $\sigma: Y \longrightarrow X$ be the blow-up of $X$ at $S=\lbrace p_1, \dots, p_n \rbrace$, where the $p_i$'s are mutually distinct points of $X$. The following result is the technical heart of the article.

\begin{thm}\label{linebundle}
		 There exists an ample line bundle $\Delta \in \Pic(X)$ such that 
        \[ C^2 \geq -\frac{1}{2} \big( \delta(X) + (\Delta.\bar{C}) \big)-n,  \]
        for all integral curves $C \subset Y$ such that $\bar{\kappa}(X \setminus \bar{C}) \geq 0$. Here, $\bar{C}:=\sigma(C)$, $\delta(X):= 3e(X) -K_X^2$ is the Miyaoka-Yau number, and $\bar{\kappa}$ denotes the logarithmic Kodaira dimension.
\end{thm}

\begin{proof}
	Let us assume that our curve $C$ is not one of the exceptional divisors. The projection of $C$ to $X$ is $\bar{C}:= \sigma(C)$. By pulling-back to $Y$, we see that $\sigma^*\bar{C} = C + \IE$, where $\IE= \sum_{i=1}^n m_iE_i$ is the total exceptional divisor coming from the multiplicities of $\bar{C}$ at the $p_i$'s.
    
	We can write the elements of $S$ as follows
    \[ S = \lbrace q_1, \dots, q_s, q_1', \dots, q_t', q_1'', \dots, q_v'' \rbrace, \]
    where $q_1, \dots , q_s \in {\rm Sing}(\bar{C})$, $q_1' , \dots, q_t' \in \bar{C} \setminus {\rm Sing}(\bar{C})$, and $q_1'', \dots, q_v'' \notin \bar{C}$. Then,
    \begin{align*}
   		C^2 &= \bar{C}^2 - \IE^2 = \bar{C}^2 - \sum_{p \in S} m_p(\bar{C})^2 \\
        &= \bar{C}^2 -  \sum_{i=1}^s m_{q_i}^{2}(\bar{C}) - \sum_{j=1}^t m_{q_j'}^{2}(\bar{C}) - \sum_{k=1}^v m_{q_k''}^{2}(\bar{C})\\
        &= \bar{C}^2 - \sum_{i=1}^s m_{q_i}^{2}(\bar{C}) - t.
    \end{align*}
    
Using Theorem \ref{SOZgen}, one gets
    \begin{align*}
    	C^2 &\geq -\frac{1}{2}\big( \delta(X) + K_X.\bar{C} \big) +\sum_{p \in {\rm Sing}(\bar{C})} \Big(1 + \frac{1}{2m_p(\bar{C})} \Big) \mu_{p}(\bar{C}) - \sum_{i=1}^s m_{q_i}^{2}(\bar{C}) - t \\        
       &=-\frac{1}{2}\big( \delta(X) + K_X.\bar{C} \big)   
       +\sum_{i=1}^s \Big(1+ \frac{1}{2m_{q_i}(\bar{C})} \Big)\mu_{q_i}(\bar{C})           
       +\sum_{p \in {\rm Sing}(\bar{C})\setminus S} \Big(1 + \frac{1} {2m_p(\bar{C})} \Big) \mu_{p}(\bar{C}) \\
       & - \sum_{i=1}^s m_{q_i}^{2}(\bar{C}) - t.
        \end{align*}
        
       Let us observe that
       \begin{align*}
       		&\sum_{i=1}^s \Big( 1+ \frac{1}{2m_{q_i}(\bar{C})} \Big)\mu_{q_i}(\bar{C}) - \sum_{i=1}^s m_{q_i}^{2}(\bar{C}) = \sum_{i=1}^s \Bigg[ \Big( 1+ \frac{1}{2m_{q_i}(\bar{C})} \Big) \mu_{q_i}(\bar{C}) - m_{q_i}^{2}(\bar{C}) \Bigg]\\
            = & \sum_{i=1}^s \Bigg[ \Big( 1+ \frac{1}{2m_{q_i}(\bar{C})} \Big) \big(\mu_{q_i}(\bar{C})-m_{q_i}(\bar{C})^2 \big) + \frac{m_{q_i}(\bar{C})}{2} \Bigg]\\
            \geq & \sum_{i=1}^s \Bigg[ \Big( 1+ \frac{1}{2m_{q_i}(\bar{C})} \Big) \big(1-2m_{q_i}(\bar{C}) \big) + \frac{m_{q_i}(\bar{C})}{2} \Bigg]\\
            = & \sum_{i=1}^s \frac{1- 3m_{q_i}(\bar{C})^2}{2m_{q_i}(\bar{C})},
    \end{align*}
where in the inequality above we have used that $\mu_p(\bar{C}) \geq \big(m_p(\bar{C})-1 \big)^2$ for every isolated singularity $p \in \bar{C}$ (see for instance \cite[Theorem 1.8]{Lin}).

From this, we deduce that
	\begin{align*}
    	C^2 &\geq -\frac{1}{2}\big( \delta(X) + K_X.\bar{C} \big)  + \sum_{i=1}^s \frac{1- 3m_{q_i}(\bar{C})^2}{2m_{q_i}(\bar{C})} - t +\sum_{p \in {\rm Sing}(\bar{C})\setminus S} \Big(1 + \frac{1}{2m_p(\bar{C})} \Big) \mu_{p}(\bar{C}) \\
        & \geq -\frac{1}{2}\big( \delta(X) + K_X.\bar{C} \big) + \sum_{i=1}^s \frac{1- 3m_{q_i}(\bar{C})^2}{2m_{q_i}(\bar{C})} - t\\ 
        &\geq -\frac{1}{2}\big( \delta(X) + K_X.\bar{C} \big)  + \sum_{i=1}^s \frac{1 - 3m_{q_i}(\bar{C})^2}{2m_{q_i}(\bar{C})} - n \quad\quad (\text{since } n\geq t)\\
        & \geq -\frac{1}{2}\big( \delta(X) + K_X.\bar{C} \big)  - \sum_{i=1}^s  \frac{3}{2}m_{q_i}(\bar{C})  - n.
    \end{align*}	
    
    At this point, we need to get rid of the multiplicities, by replacing them with suitable intersection numbers. Let us choose a very ample line bundle $A \in \Pic(X)$, and let $\varphi_A: X \longra \IP^{h^0(A)-1}$ be the corresponding embedding. Then, the multiplicities $m_{p_i}(\bar{C})$ are bounded by the degree of $\bar{C}$ in the embedding $\varphi_A$, i.e.~$m_{p_i}(\bar{C}) \leq (\bar{C}.A)=\deg_{\IP^{h^0(A)-1}}(\bar{C})$. Therefore, it follows that
    
    	\begin{align*}
    	C^2 &\geq -\frac{1}{2}\big( \delta(X) + K_X.\bar{C} \big)  - \sum_{i=1}^s  \frac{3}{2}m_{q_i}(\bar{C})  - n\\
        &\geq -\frac{1}{2}\big( \delta(X) + K_X.\bar{C} \big)  - \sum_{i=1}^s  \frac{3}{2}(\bar{C}.A) - n \quad\quad (\text{since } n\geq s)\\
        &\geq -\frac{1}{2} \big( \delta(X) + (K_X+3nA).\bar{C} \big) -n.
    \end{align*}
    
    The line bundle $K_X+3nA$ might not be ample, but it becomes such upon replacing $A$ with a multiple. This means that for a suitable choice of $A$, the adjoint line bundle $\Delta:= K_X + 3nA$ is ample, thus
    \[ C^2 \geq - \frac{1}{2} \big( \delta(X) + \Delta.\bar{C} \big) - n \]
    
    This concludes the proof in case $C$ is not one of the exceptional divisor. However, if $C$ were to be one of the exceptional divisors, the bound above would still hold true, therefore we are done.
\end{proof}

As a consequence, we immediately get a linear bound on the self-intersection of integral curves on all surfaces $Y$ as above having the additional requirement that their Kodaira dimension is non-negative.

\begin{cor}\label{thmA}
		Assume $X$ is a surface of non-negative Kodaira dimension. Then, in the setting above, there exists a big and nef line bundle $\Gamma$ that bounds negativity, i.e.,
        \[ C^2 \geq -\frac{1}{2} \big( \delta(X) + C.\Gamma \big) - n, \]
        for every integral curve $C \subset Y$. In other words, if we define $\deg_\Gamma C:= (C.\Gamma)$, then
  \[ C^2 \geq -\Bigg( \frac{1}{2}\delta(X) + n \Bigg) - \frac{1}{2}\deg_\Gamma C, \]
  i.e.~the negativity of $C$ is bounded by a function that depends on $X$, the number of points we have blown-up, and the $\Gamma$-degree of $C$.
\end{cor} 

\begin{proof}
The line bundle $\Delta$ in the proof of Theorem \ref{linebundle} provides us with a degree function on $\NS(X)$. As a consequence, we obtain a choice of a degree-like line bundle of $Y$ by setting $\Gamma:= \sigma^*\Delta$. The line bundle $\Gamma$ will never be ample (we are pulling back along a blow-up), but it is nevertheless big and nef. Hence we can use it to provide a weighted bound for the negativity on $Y$.
\end{proof}

It is interesting to observe the following facts:

\begin{itemize}
	\item if $X$ is a minimal surface, then the bound of the negativity of $Y$ directly arises naturally from its minimal model;
    
    \item the bound on the negativity is now linear in $(C.\Gamma)$, while the weighted BNC predicts the existence of a quadratic bound.
\end{itemize} 

\section{Bounding negativity on blow-ups of $\IP^2$}

In this section, we will study the problem of bounding negativity for blow-ups of $\IP^2$. We present here two different approaches to find bounds for the intersection numbers for curves on blow-ups of the complex projective plane. We start with the first approach using Orevkov-Sakai-Zaidenberg's inequality.

\begin{thm}\label{boundrat}
		Let $\sigma: Y \longra \IP^2$ be the blow-up of $\IP^2$ at $S= \lbrace p_1, \dots, p_n \rbrace$, where the $p_i$'s are distinct points of $\IP^2$, and let $C$ be an irreducible and reduced curve on $Y$. Then,
        \[ C^2 \geq - 2n(C.L), \]
        where $L$ is the pull-back of a line in $\IP^2$.
\end{thm}

\begin{proof}
In this case, there do exist curves for which the logarithmic Kodaira dimension of the complement is $-\infty$. As it was shown by Wakabayashi \cite{Waka}, if $D \subset \mathbb{P}^{2}$ is an irreducible and reduced curve of degree $d \geq 4$ having $s\geq 1$ singular points, which is not a rational cuspidal curve with one cusp, then the logarithmic Kodaira dimension of $\mathbb{P}^{2} \setminus D$ is non-negative. Therefore, we can apply Theorem \ref{SOZgen} to bound the self-intersection of these curves. In fact, it was pointed by Sakai \cite{Sakai} that the inequality in Theorem \ref{SOZ} holds for \emph{all} irreducible and reduced curves $D \subset \mathbb{P}^2$ of degree $d\geq 3$ -- it is enough to verify the remaining cases by simple computations. 

Let $C \subset Y$ be an irreducible an reduced curve, and let us denote by $\bar{C}$ its image under $\sigma$. If $\bar{C}.H \geq 3$, $H$ being the class of a line in $\IP^2$, then we can repeat the proof of Theorem \ref{linebundle} to obtain
\begin{align*}
		C^2 &\geq -\frac{1}{2}\big( \delta(\IP^2) + K_{\IP^2}.\bar{C} \big)  -  \frac{3}{2}\sum_{i=1}^s m_{q_i}(\bar{C})  - n\\
        & = \frac{3}{2}(H.\bar{C}) - \frac{3}{2}\sum_{i=1}^s m_{q_i}(\bar{C})  - n\\
        &\geq \frac{3}{2}(H.\bar{C}) -  \frac{3}{2}\sum_{i=1}^s (H.\bar{C}) - n\\
        &\geq -\frac{3}{2}n(H.\bar{C})-n \geq -2n(H.\bar{C}) = -2n(L.C),
\end{align*}
where $L=\sigma^*H$. We are left to deal with curves $C \subset Y$ whose image $\bar{C}$ is either a line or a conic. For such curves, we have that 
\[ 1 - \sum_{i=1}^s m_{p_i}(\bar{C}) \leq C^2 \leq 2 - \sum_{i=1}^s m_{p_i}(\bar{C}). \]
However, due to the restriction on the degree, $\bar{C}$ is necessarily smooth and $m_{p_i}(\bar{C}) =1$ for all $i=1, \dots, s$. Therefore, we find that $C^2 \geq 1-n$, and thus we have proven the result.
\end{proof}

Our second approach to the problem allows us to improve our previous bound from Theorem \ref{boundrat}, and this is a consequence of a classical result in the theory of algebraic curves \cite[Theorem 7.22]{Wall}. 
\begin{thm}(Pl\"ucker-Teissier formula)
Let $C \subset \mathbb{P}^{2}$ be an irreducible and reduced curve. Then
$$\sum_{p \in {\rm Sing}(C)}(\mu_{p} + m_{p} - 1) \leq d(d-1).$$
\end{thm}

In the setting of Theorem \ref{boundrat}, by using the inequality $\mu_p \geq (m_p(\bar{C})-1)^2$ for $p \in \text{Sing}(\bar{C})$, the Pl\"ucker-Teissier formula implies that (again, we use the notation as in the proof of Theorem \ref{linebundle}): 
\[ d^2 - d \geq \sum_{p \in {\rm Sing}(\bar{C})} (\mu_p(\bar{C})+m_p(\bar{C}) -1) \geq \sum_{i=1}^s (\mu_{q_i}(\bar{C})+m_{q_i}(\bar{C}) -1) \geq \sum_{i=1}^s m_{q_i}(\bar{C})\big( m_{q_i}(\bar{C})-1 \big),\]
which in turn shows that 
\begin{align*}
C^2 &= d^2 - \sum_{p \in S} m_p(\bar{C})^2 = d^2 - \sum_{i=1}^s m_{q_i}(\bar{C})^2 -t\\
&\geq d - \sum_{i=1}^s m_{q_i}(\bar{C}) -t \geq d(1-s) -t \geq -d(s+t) \geq -nd,
\end{align*}
and we got a better constant than in the statement of Theorem \ref{boundrat}.\\

We would like to conclude by making the following remark, which considers the case of a blow-up of $\IP^2$ at a set $\mathcal{P}$ of points in very general position. Assume that $\mathcal{P} = \{p_{1}, ...,p_{n}\}$ are points in very general position and we consider the blowing-up $\pi : X \rightarrow \mathbb{P}^{2}$ along $\mathcal{P}$. Let $C \subset X$ be an irreducible and reduced curve, and denote by $\bar{C}\subset\mathbb{P}^{2}$ its image. Then by \cite[Lemma 1]{Xu}, one has:
$${C}^{2} \geq - \min \{m_{q_{1}}, \dots , m_{q_{s}},m_{q'_1}, \dots , m_{q'_t} \} \geq -d,$$
which means that in generic case the better bound $C^2 \geq -d$ holds for every irreducible and reduced curve $C \subset X$. Notice that this bound does not depend on the number of points that we have blown up the surface.

\section{Bounding negativity on blow-ups of Hirzebruch surfaces}\label{hirzebruch}

We denote by $\IF_m$ the $m$\textsuperscript{th} Hirzebruch surface, and let us consider the case $m \neq 1$ only, so that $\IF_m$ is a minimal surface ($\IF_1$ is $\IP^2$ blown-up at one point). If $F$ is the class of a fiber, and $H$ is the tautological section of $\IF_m$, then 
\[ \Pic(\IF_m) = \Z F \oplus \Z H, \qquad H^2 = m, \qquad H.F=1, \qquad K_{\IF_m} = -2H+(m-2)F. \]
We would like to mimic the argument for blow-ups of $\IP^2$. Let $\sigma: Y \longra \IF_m$ be the blow-up of $\IF_m$ at a set $S = \lbrace p_1, \dots, p_n \rbrace$ of distinct points. Suppose that $C \subset Y$ is a curve with the property that $\bar{\kappa}(Y \setminus C) \geq 0$, and let $\bar{C}$ be its image under $\sigma$. By the proof of Theorem \ref{linebundle}, we get
\[ C^2 \geq -\frac{1}{2}\big( -4 + K_X.\bar{C} \big) -\frac{3}{2} \sum_{i=1}^s m_{p_i}(\bar{C}) - n. \]

Now, the line bundle $A:=H+F$ is very ample by \cite[Exercise IV.18(2)]{Beauville}, and it embeds $\IF_m$ into $\IP^{m+3}$ as a surface of degree $m+2$. Therefore,
\begin{align*}
		C^2 &\geq 2-n + \frac{1}{2} \Big( -K_X.\bar{C} - 3 \sum_{i=1}^s m_{p_i}(\bar{C}) \Big)\\
        &\geq 2-n - \frac{1}{2} \Big( K_X.\bar{C} + 3 \sum_{i=1}^s A.\bar{C} \Big)\\
        &\geq 2-n - \frac{1}{2} \Big( K_X.\bar{C} + 3n (A.\bar{C}) \Big)\\
        &\geq 2-n - \frac{1}{2}\big((K_X+3nA).\bar{C}\big).
\end{align*}

The line bundle $\Delta:= K_X + 3nA$ is always very ample on $\IF_m$, thus yielding a big and nef line bundle $\Gamma := \sigma^*\Delta$ on $Y$ that bounds the negativity on $Y$:
\[ C^2 \geq 2-n - \frac{1}{2}\deg_\Gamma C.\]

It is natural to ask for which classes of curves we can apply our lower-bound, and the answer is provided by the following Wakabayashi-type result \cite[Theorem 1.4]{Moe}.
\begin{thm}
On a Hirzebruch surface $\IF_m$, let $C$ be an irreducible curve of genus $g$ and type $(a,b)$ with $b > 2$, $a > 2 - \frac{1}{2}bm$, and $a\geq 0$. Then
\begin{itemize}
\item If $g > 0$, then the logarithmic Kodaira dimension of $\IF_m \setminus C$ is equal to $2$.
\item If $g=0$ and $C$ has at least three cusps, then the logarithmic Kodaira dimension of $\IF_m \setminus C$ is equal to $2$.
\item If $g=0$ and $C$ at least two cusps, then the logarithmic Kodaira dimension of $\IF_m \setminus C$ is at least equal to $0$.
\end{itemize}
\end{thm}

\vspace{5mm} 
\paragraph*{\bf{Acknowledgement.}}
The idea behind this paper grew up during the Workshop \emph{Newton-Okounkov Bodies, Test Configurations, and Diophantine Geometry} in Banff in 2017. We would like to thank Alex K\"uronya, Mike Roth, and Tomasz Szemberg for organizing an excellent workshop, and Joaquim Ro\'e and Alex K\"uronya for fruitful discussions. We also would like to thank Xavier Roulleau and Mikhail Zaidenberg for comments. We would like warmly thank \emph{Banff International Research Center} for the extraordinary hospitality. Finally, we would like to warmly thank an anonymous referee for valuable comments. The second author was partially supported by the Fundation for Polish Science (FNP) Scholarship Start No. 076/2018.



\end{document}

%% file: pack.tex
\usepackage{amsmath}
\usepackage{amsfonts}
\usepackage{amssymb}
\usepackage{amsthm}
\usepackage{mathtools}
\usepackage{mathrsfs}
\usepackage{float}
\usepackage[all,cmtip]{xy}
      \theoremstyle{definition}
\newtheorem{defi}{Definition}[section]
\theoremstyle{plain}
\newtheorem{thm}[defi]{Theorem}

\newtheorem{cor}[defi]{Corollary}

\theoremstyle{remark}
\newtheorem{constr}[defi]{Construction}

\newtheorem{conj}[defi]{Conjecture}

\theoremstyle{definition}

\newtheorem*{theoremA}{Theorem A}
\newtheorem*{theoremB}{Theorem B}

\newcommand{\longra}{\longrightarrow}

\newcommand{\IE}{{\mathbb E}}
\newcommand{\IF}{{\mathbb F}}

\newcommand{\IP}{{\mathbb P}}
\newcommand{\Q}{{\mathbb Q}}

\newcommand{\Z}{{\mathbb Z}}

\newcommand{\Pic}{\operatorname{Pic}}

\newcommand{\NS}{\operatorname{NS}}

  \makeatletter
\newcommand{\xdashrightarrow}[2][]{\ext@arrow 0359\rightarrowfill@@{#1}{#2}}
\newcommand{\xdashleftarrow}[2][]{\ext@arrow 3095\leftarrowfill@@{#1}{#2}}
\newcommand{\xdashleftrightarrow}[2][]{\ext@arrow 3359\leftrightarrowfill@@{#1}{#2}}
\def\rightarrowfill@@{\arrowfill@@\relax\relbar\rightarrow}
\def\leftarrowfill@@{\arrowfill@@\leftarrow\relbar\relax}
\def\leftrightarrowfill@@{\arrowfill@@\leftarrow\relbar\rightarrow}
\def\arrowfill@@#1#2#3#4{%
  $\m@th\thickmuskip0mu\medmuskip\thickmuskip\thinmuskip\thickmuskip
   \relax#4#1
   \xleaders\hbox{$#4#2$}\hfill
   #3$%
}
\makeatother
